\documentclass[12pt]{amsart}
\usepackage{amsthm, amssymb}
\usepackage{amsfonts, amscd}
\usepackage{epsfig,multicol}
\usepackage{xypic}
\theoremstyle{plain} \numberwithin{equation}{section}
\newtheorem{thm}{Theorem}[section]
\newtheorem{cor}[thm]{Corollary}
\newtheorem{conj}[thm]{Conjecture}
\newtheorem{prob}[thm]{Problem}
\newtheorem{prop}[thm]{Proposition}

\theoremstyle{definition}
\newtheorem{defn}{Definition}

\theoremstyle{remark}
 \newtheorem{rem}{Remark}

\def\Z{\Bbb Z}

\begin{document}
\title[Equivariant bordism of 2-torus manifolds and unitary toric manifolds]{\large \bf Equivariant bordism of 2-torus manifolds and unitary toric manifolds--a survey}
\author[Zhi L\"u]{Zhi L\"u }
\footnote[0]{
\endgraf
 {\bf 2000AMS Classification:} 57S10, 57R85, 14M25, 52B70.
 \endgraf
 Supported by grants from NSFC (No. 11371093, No. 11431009 and No. 11661131004).}
\address{School of Mathematical Sciences, Fudan University, Shanghai,
200433, P.R. China.} \email{zlu@fudan.edu.cn}


\begin{abstract}
In this paper we survey results and  recent progresses on the equivariant bordism classification of 2-torus manifolds and unitary toric manifolds.
\end{abstract}

\maketitle


\section{Introduction}\label{int}
 An $n$-dimensional {\em 2-torus
manifold} is a smooth closed $n$-dimensional (not necessarily
oriented) manifold equipped with an  effective smooth  $({\Bbb Z}_2)^n$-action,
so its  fixed point set  is empty or consists
 of a set of isolated points (see \cite{l1, lm}). A $2n$-dimensional  {\em  unitary toric manifold},  introduced by Masuda in~\cite{ma},  is a unitary $2n$-dimensional manifold equipped with an effective $T^n$-action fixing a nonempty fixed point set and preserving the tangential stably complex structure of $M^{2n}$,  where $T^n$ is the torus group of rank $n$ and a unitary manifold is an   oriented closed smooth manifold  whose tangent bundle admits a stably complex structure.
The seminal work of Davis and Januszkiewicz in
\cite{dj} studied two kinds of equivariant manifolds: {\em small covers}  and {\em quasitoric manifolds}, which are   the  real and complex topological versions of  toric
varieties in algebraic geometry, respectively, where a small cover of dimension $n$ (resp. a quasitoric manifold of dimension $2n$) is a smooth closed $n$-dimensional manifold (resp. $2n$-dimensional manifold)  admiting a locally standard $({\Bbb Z}_2)^n$-action (resp. $T^n$-action)  such that the orbit space  is  a simple convex polytope. As shown in~\cite{dj}, small covers and quasitoric manifolds have a very beautiful algebraic topology and provide a strong link between equivariant topology, polytope theory and combinatorics. Obviously, each small cover is a special 2-torus manifold. Buchstaber and Ray showed in~\cite{br}  that each quasitoric manifold with an omniorientation always admits a compatible tangential stably complex structure.  Thus, small covers and omnioriented quasitoric manifolds provide abundant examples of 2-torus manifolds and unitary toric manifolds, respectively.

\vskip .1cm
In nonequivariant case, Buchstaber and Ray showed in~\cite{br} that each class of $\mathfrak{N}_n$ (resp. $\Omega_{2n}^U$) contains an $n$-dimensional  small cover (resp. a $2n$-dimensional quasitoric manifold) as its representative, where $\mathfrak{N}_*=\sum_{m\geq 0} \mathfrak{N}_m$ (resp. $\Omega_*^U=\sum_{m\geq 0}\Omega_{2m}^U$) is the ring formed by the unoriented bordism classes of all smooth closed manifolds (resp. the unitary bordism classes of all unitary manifolds). The work of Buchstaber and Ray gives a motivation to the study of the equivariant bordism classification of 2-torus manifolds and unitary toric manifolds, so that  the following question arises naturally.
\begin{enumerate}
\item[(Q1)] {\em Can preferred representatives in the
classes of $\mathcal{Z}_n(({\Bbb Z}_2)^n)\ ($resp. $\mathcal{Z}_{2n}^{U}(T^n))$ be chosen from small covers (resp. omnioriented quasitoric manifolds)?}
\end{enumerate}
 where $\mathcal{Z}_n(({\Bbb Z}_2)^n)\ ($resp. $\mathcal{Z}_{2n}^{U}(T^n))$ denotes the group produced by the $({\Bbb Z}_2)^n$-equivariant unoriented bordism classes of all $n$-dimensional 2-torus manifolds (resp. the $T^n$-equivariant unitary bordism classes of all $2n$-dimensional unitary toric manifolds). Note that the cartesian products of actions also define the graded rings $\mathfrak{M}_*=\bigoplus_{n\geq 0}\mathcal{Z}_n(({\Bbb Z}_2)^n)$ and $\Xi_*=\bigoplus_{n\geq 0}\mathcal{Z}_{2n}^{U}(T^n)$, which, significantly, turn out to be non-commutative.
  With respect to this question, the author of this paper first dealt with the case of 2-torus manifolds and proposed the following concrete conjecture  in~\cite{l1}.
\vskip .1cm
\noindent {\bf Conjecture $(\ast)$:} {\em Each class of $\mathcal{Z}_n(({\Bbb Z}_2)^n)$ contains a  small cover  as its representative.}
\vskip .1cm
\noindent It was shown in~\cite{l1} that Conjecture $(\ast)$ is true for $n\leq 3$, and there is also an essential link among the equivariant bordism of 3-dimensional 2-torus manifolds, 3-dimensional colored polytopes and mod 2 GKM graphs of valence 3. Moreover, the following question also arises naturally although 2-torus manifolds (resp. unitary toric manifolds) form a much wider class than small covers (resp. omnioriented quasitoric manifolds).
\begin{enumerate}
\item[(Q2)] {\em Is there still a strong link among the equivariant bordism of 2-torus manifolds and unitary toric manifolds, colored polytopes and (mod 2) GKM graphs in the general case?}
\end{enumerate}

The purpose of this paper is to investigate results and  recent development on the equivariant bordism classification of 2-torus manifolds and unitary toric manifolds with respect to the above questions.

\vskip .1cm
 In the setting of 2-torus manifolds, a  significant amount of technical machinery has been developed in~\cite{lt}  by defining  a differential operator on the ¡°dual¡± algebra of the unoriented $({\Bbb Z}_2)^n$-representation algebra introduced by Conner and Floyd, so that the satisfactory solutions  to the  questions (Q1) and (Q2) can be obtained, and in particular, Conjecture $(*)$ can be answered affirmatively.  In addition, this technical machinery can be combined with the mod 2 GKM theory and the Davis--Januszkiewicz Theory of small covers together very well, so that one can determine how the
graded noncommutative ring $\mathfrak{M}_*$  is generated, and find some essential relationships among 2-torus
manifolds, coloring polynomials, colored simple convex polytopes, colored graphs.  It should be pointed out that the classical  $({\Bbb Z}_2)^n$-equivariant bordism theory and results (e.g., tom Dieck's existence theorem)  also play important roles on the study of the equivariant bordism classification of 2-torus manifolds. In Section~\ref{2-torus}, we shall systemically introduce the developed equivariant unoriented bordism theory in the setting of 2-torus manifolds.
\vskip .1cm

In the setting of unitary toric manifolds,  whether or not  each class of $\mathcal{Z}_{2n}^{U}(T^n)$ is represented by an omnioriented quasitoric manifold is still open in the general case. However, some techniques and ideas developed in the setting of 2-torus manifolds can still be carried out very well in this case. Indeed, Darby in~\cite{dar} drew those techniques and ideas  into the case of unitary toric manifolds with some new viewpoints, so that significant advances can be obtained in many aspects, such as  the graded noncommutative ring $\Xi_*$, some essential links to polytope theory and torus graphs, and so on. The key reason why the case of unitary toric manifolds has  less progresses than that of 2-torus manifolds is because of 
the existence of an  infinite number of irreducible complex $T^n$-representations. This leads the case of unitary toric manifolds to be much more difficult and complicated. We shall give an introduction on the study of $\mathcal{Z}_{2n}^{U}(T^n)$ in Section~\ref{unitary}. In addition, in Section~\ref{unitary} we also shall mention the work on the equivariant Chern numbers of unitary toric manifolds in~\cite{lt1} which is  related to Kosniowski conjecture, and some problems and conjectures will also be proposed therein.
\vskip .1cm
Finally we survey a result on the relation between $\mathcal{Z}_n(({\Bbb Z}_2)^n)$ and $\mathcal{Z}_{2n}^{U}(T^n))$ in Section~\ref{rela}.

\section{Equivariant unoriented bordism  of 2-torus manifolds}\label{2-torus}

 In the early 1960s, Conner and Floyd~\cite{cf} (also see \cite{co}) begun the study of geometric equivariant unoriented and oriented bordism theories for  smooth closed manifolds with periodic diffeomorphisms, and the  subject has also continued to develop and to flourish by extending their ideas to other equivariant bordisms  since then. For example, the homotopy theoretic analogue was described by tom Dieck~\cite{tom3}.  A main aspect of Conner--Floyd's work is on the study of $({\Bbb Z}_2)^n$-equivariant unoriented bordism. Conner and Floyd studied the localization of $({\Bbb Z}_2)^n$-equivariant unoriented bordism. The following localization theorem is due to Conner and Folyd~\cite{cf} for $n=1$ and Stong~\cite{s} for the general case.
\begin{thm} \label{ring}
The ring homomorphism
$\phi_*: \mathfrak{N}_*^{({\Bbb Z}_2)^n}=\sum_{i\geq 0}\mathfrak{N}_i^{({\Bbb Z}_2)^n}\longrightarrow \mathfrak{N}_*(BO),$
defined by mapping the equivariant unoriented class of a smooth closed $({\Bbb Z}_2)^n$-manifold to the unoriented bordism class of its normal bundle of fixed point set,  is injective, where $\mathfrak{N}_*^{({\Bbb Z}_2)^n}$ is the ring formed by the equivariant unoriented bordism classes of all smooth closed $({\Bbb Z}_2)^n$-manifolds.
\end{thm}
However, generally speaking,   the ring structure of $\mathfrak{N}_*^{({\Bbb Z}_2)^n}$ is still far from settled except for the case $n=1$ (see~\cite{a, sinh}). In~\cite{cf}, Conner and Floyd discussed the case in which the fixed point set of the action is isolated, and introduced and studied a graded commutative algebra  over ${\Bbb
Z}_2$  with unit, $\mathcal{Z}_*(({\Bbb Z}_2)^n)=\sum_{m\geq 0}\mathcal{Z}_m(({\Bbb Z}_2)^n)$,
where $\mathcal{Z}_m(({\Bbb Z}_2)^n)$ consists of $({\Bbb Z}_2)^n$-equivariant unoriented
bordism  classes of all smooth closed $m$-manifolds with effective  $({\Bbb Z}_2)^n$-actions fixing a finite set. Clearly,  $\mathcal{Z}_*(({\Bbb Z}_2)^n)$ is a subring of $\mathfrak{N}_*^{({\Bbb Z}_2)^n}$, and  when $m=n$,  $\mathcal{Z}_n(({\Bbb Z}_2)^n)$ is exactly formed by the classes of all 2-torus $n$-manifolds. Also, the restriction to $\mathcal{Z}_*(({\Bbb Z}_2)^n)$ of $\phi_*$  gives the following monomorphism (still denoted by $\phi_*$)
$$\phi_*: \mathcal{Z}_*(({\Bbb Z}_2)^n)\longrightarrow \mathcal{R}_*(({\Bbb Z}_2)^n)$$ defined
 by $\{M\}\longmapsto \sum_{p\in M^{({\Bbb Z}_2)^n}}[\tau_pM]$ where $\tau_pM$ denotes
the real $({\Bbb Z}_2)^n$-representation on the tangent space at $p\in
M^{({\Bbb Z}_2)^n}$, and $\mathcal{R}_*(({\Bbb Z}_2)^n)=\sum_{m\geq 0}\mathcal{R}_m(({\Bbb Z}_2)^n)$  is  the
  graded polynomial algebra over ${\Bbb Z}_2$ generated by the  isomorphism classes of one-dimensional irreducible  real $({\Bbb Z}_2)^n$-representations, which was introduced by Conner and Folyd and is called the {\em Conner--Floyd unoriented $({\Bbb Z}_2)^n$-representation algebra} here.  Conner and Floyd showed in \cite{cf} that
when $n=1$ $\mathcal{Z}_*(\Z_2)\cong\Z_2$, and when $n=2$, $\mathcal{Z}_*((\Z_2)^2)\cong\Z_2[u]$ where $u$ denotes the class of
${\Bbb R}P^2$ with the standard  $(\Z_2)^2$-action. Since then, any new progress on $\mathcal{Z}_*(({\Bbb Z}_2)^n)$ has not been made until the work in ~\cite{l1} appeared in 2009. When
$n=3$, the group structure of $\mathcal{Z}_3((\Z_2)^3)$ was determined in \cite{l1}
(see also \cite{ly}), and it was also shown therein that $\dim_{\Z_2}\mathcal{Z}_3((\Z_2)^3)=13$.

\vskip .1cm The objective of this section is to survey the recent progress on the study of $\mathcal{Z}_n(({\Bbb Z}_2)^n)$ (i.e., the equivariant unoriented bordism of 2-torus manifolds) in the general case.

\subsection{The reformulation of the existence theorem of tom Dieck and a differential operator}
In \cite[Theorem 6]{tom1}, tom Dieck showed an existence theorem, saying that  the existence of an $m$-dimensional smooth closed $({\Bbb Z}_2)^n$-manifold $M^m$
fixing a finite set can be characterized by the integral property of its fixed point data.  In~\cite{lt}, L\"u and Tan gave a simple proof to show that the
existence theorem of tom Dieck can be  formulated  into  the
following result in terms of Kosniowski and Stong's
localization formula (\cite{ks}).
\begin{thm} [{\cite[Theorem 2.2]{lt}}] \label{dks}
Let $\{\tau_1, ..., \tau_l\}$ be a collection of $m$-dimensional faithful
$({\Bbb Z}_2)^n$-representations in $\mathcal{R}_m(({\Bbb Z}_2)^n)$. Then a necessary and
sufficient condition that $\tau_1+\cdots+\tau_l\in \text{\rm Im}\phi_m$ $($or
$\{\tau_1, ..., \tau_l\}$ is the fixed point data of a
$({\Bbb Z}_2)^n$-manifold $M^m)$ is that for all symmetric polynomial functions
$f(x_1,...,x_m)$ over ${\Bbb Z}_2$,
\begin{equation} \label{formula-tks1}\sum_{i=1}^l{{f(\tau_i)}\over{\chi^{({\Bbb Z}_2)^n}(\tau_i)}}\in
H^*(B({\Bbb Z}_2)^n;{\Bbb Z}_2)\end{equation}  where $\chi^{({\Bbb Z}_2)^n}(\tau_i)$
denotes the equivariant Euler class of $\tau_i$, which is a product
of $m$ nonzero elements of $H^1(B({\Bbb Z}_2)^n;{\Bbb Z}_2)$, and $f(\tau_i)$
means that variables $x_1,...,x_m$ in the function $f(x_1,...,x_m)$
are replaced by those $m$ degree-one factors in
$\chi^{({\Bbb Z}_2)^n}(\tau_i)$.
\end{thm}
\begin{rem}
  Although all elements of $\text{\rm Im} \phi_*$ can be characterized by the formula (\ref{formula-tks1}),
  it is still quite difficult to determine the algebra structure of $\text{\rm Im}
\phi_*\cong \mathcal{Z}_*(({\Bbb Z}_2)^n)$.
\end{rem}

It is well-known that all irreducible
real $({\Bbb Z}_2)^n$-representations bijectively correspond to all  elements
in $\text{\rm Hom}(({\Bbb Z}_2)^n,{\Bbb Z}_2)$, where every irreducible real representation
of $({\Bbb Z}_2)^n$ has the form $\lambda_\rho: ({\Bbb Z}_2)^n\times{\Bbb
R}\longrightarrow{\Bbb R}$ with $\lambda_\rho(g,x)=(-1)^{\rho(g)}x$
for $\rho\in\text{\rm Hom}(({\Bbb Z}_2)^n,{\Bbb Z}_2)$, and $\lambda_\rho$ is trivial if
$\rho(g)=0$ for all $g\in ({\Bbb Z}_2)^n$. Write $J_n^{\Bbb R}=\text{\rm Hom}(({\Bbb Z}_2)^n,{\Bbb Z}_2)$ and regard it as the set of all all irreducible real $({\Bbb Z}_2)^n$-representations. Then the free polynomial algebra on $J_n^{\Bbb R}$ over ${\Bbb Z}_2$, denoted by ${\Bbb Z}_2[J_n^{\Bbb R}]$, can be identified with
 $\mathcal{R}_*(({\Bbb Z}_2)^n)$.
 Similarly, one has also another free polynomial algebra ${\Bbb Z}_2[J_n^{*{\Bbb R}}]$ on $J_n^{*{\Bbb R}}$ over ${\Bbb Z}_2$, which is called the {\em dual algebra} of ${\Bbb Z}_2[J_n^{\Bbb R}]=\mathcal{R}_*(({\Bbb Z}_2)^n)$, where $J_n^{*{\Bbb R}}=\text{\rm Hom}({\Bbb Z}_2, ({\Bbb Z}_2)^n)$ is the dual of $J_n^{\Bbb R}$ as ${\Bbb Z}_2$-linear spaces.
 \begin{defn}
 A square-free homogeneous polynomial  $g =\sum_it_{i,1}\cdots t_{i,n}$   of degree $n$ in ${\Bbb Z}_2[J_n^{\Bbb R}]$ is called a {\em faithful $({\Bbb Z}_2)^n$-polynomial} if each monomial $t_{i,1}\cdots t_{i,n}$ is a faithful $({\Bbb Z}_2)^n$-representative in ${\Bbb Z}_2[J_n^{\Bbb R}]$ (i.e., $\{t_{i,1}\cdots t_{i,n}\}$ forms a basis of $J_n^{\Bbb R}$).
  \end{defn}
  Both $J_n^{\Bbb R}$ and $J_n^{*{\Bbb R}}$ are isomorphic to $({\Bbb Z}_2)^n$ and are dual to each other by the following pairing:
\begin{equation}\label{pairing}
\langle\cdot, \cdot\rangle: J_n^{*{\Bbb R}}\times J_n^{\Bbb R}\longrightarrow \text{\rm Hom}({\Bbb Z}_2, {\Bbb Z}_2)
\end{equation}
defined by $\langle \xi, \rho\rangle=\rho\circ \xi$, composition of
homomorphisms. Thus, each faithful $({\Bbb Z}_2)^n$-polynomial $g =\sum_it_{i,1}\cdots t_{i,n}$   of degree $n$ in ${\Bbb Z}_2[J_n^{\Bbb R}]$
determines a unique  homogeneous polynomial $g^*=\sum_i s_{i,1}\cdots s_{i,n}$ in ${\Bbb Z}_2[J_n^{*{\Bbb R}}]$, which is called the  {\em dual
$({\Bbb Z}_2)^n$-polynomial} of $g$, where  $\{s_{i,1},..., s_{i,n}\}$ is the  dual basis of $\{t_{i,1}, ...,  t_{i,n}\}$, determined by the pairing (\ref{pairing}).

\vskip .1cm
In~\cite{lt}, L\"u and Tan defined a {\em differential operator} $d$ on ${\Bbb Z}_2[J_n^{*{\Bbb R}}]$ as follows: for each monomial $s_1\cdots s_i$ of degree
$i\geq 1$
$$d_i(s_1\cdots s_i)=\begin{cases}
\sum_{j=1}^is_1\cdots s_{j-1}\widehat{s}_js_{j+1}\cdots s_i &\text{ if } i>1\\
1 &\text{ if } i=1.
\end{cases}$$
and $d_0(1)=0$, where the symbol $\widehat{s}_j$ means that $s_j$ is
deleted. Obviously, $d^2=0$, so $({\Bbb Z}_2[J_n^{*{\Bbb R}}],
d)$ forms a chain complex. Then,  L\"u and Tan proved the following theorem which formulates a simple criterion of a faithful $({\Bbb Z}_2)^n$-polynomial $g\in \text{\rm Im}\phi_n$ in terms of the vanishing of the differential $d$ on the dual of $g$.
\begin{thm}[{\cite[Theorem 2.3]{lt}}]\label{main result}
Let $g=\sum_i t_{i,1}\cdots t_{i, n}$ be a faithful $({\Bbb Z}_2)^n$-polynomial
in ${\Bbb Z}_2[J_n^{\Bbb R}]$. Then
 $g\in \text{\rm Im}\phi_n$ if and only if $d(g^*)=0$.
\end{thm}
By $\mathcal{S}(J_n^{\Bbb R})$ (resp. $\mathcal{S}(J_n^{*{\Bbb R}})$) one denotes  the infinite symmetric tensor
algebra on $J_n^{\Bbb R}$ (resp. $J_n^{*{\Bbb R}})$ over ${\Bbb Z}_2$. It turns out that $\mathcal{S}(J_n^{\Bbb R})$ (resp. $\mathcal{S}(J_n^{*{\Bbb R}})$) is in effect the same as the graded polynomial algebra over ${\Bbb Z}_2$, in indeterminates that are basis elements for $J_n^{\Bbb R}$ (resp. $J_n^{*{\Bbb R}})$.
 Then one has that  $$\mathcal{S}(J_n^{\Bbb R})\cong \mathcal{S}(J_n^{*{\Bbb R}})
\cong H^*(B({\Bbb Z}_2)^n;\Z_2)$$ as algebras since both $J_n^{\Bbb R}$ and $J_n^{*{\Bbb R}}$ are isomorphic to $H^1(B({\Bbb Z}_2)^n;{\Bbb Z}_2)$ as ${\Bbb Z}_2$-linear spaces.
\vskip .1cm
Theorems~\ref{dks} and \ref{main result} give  following
interesting algebraic corollary, which indicates a unification of ``differential" and ``integral" in some sense.
\begin{cor}[\cite{lt}]\label{diff-formula}
Let $g=\sum_i t_{i,1}\cdots t_{i, n}$ be a faithful $({\Bbb Z}_2)^n$-polynomial
in ${\Bbb Z}_2[J_n^{\Bbb R}]$.
 Then $d(g^*)=0$ if and only if for all symmetric
polynomial functions $f(x_1,...,x_n)$ over ${\Bbb Z}_2$,
\begin{equation*}\label{integral}
\sum_{i}{{f(t_{i,1}, ..., t_{i, n})}\over{t_{i,1}\cdots t_{i,
n}}}\in\mathcal{S}(J_n^{\Bbb R})\end{equation*} when $t_{i,1}\cdots
t_{i, n}$ and $f(t_{i,1}, ..., t_{i, n})$ are regarded as
polynomials in $\mathcal{S}(J_n^{\Bbb R})$.
\end{cor}

The proof of Theorem~\ref{main result} is based upon another characterization of $g\in \text{\rm Im} \phi_n$ in terms of $({\Bbb Z}_2)^n$-colored graphs (or mod 2 GKM graphs), which will be introduced in the next subsection.

\subsection{$({\Bbb Z}_2)^n$-colored graphs and small covers}

In \cite{gkm}, Goresky, Kottwitz and MacPherson established the GKM
theory, indicating that  there is an essential link between topology
and geometry of torus actions and the combinatorics of colored
graphs (see also \cite{gz}). Such a link has already been expanded
to the case of mod 2-torus actions (see, e.g.,
\cite{bl, bgh, l2, l3}). Specifically,
assume that $M^m$ is a smooth closed $m$-manifold with an effective
smooth $({\Bbb Z}_2)^n$-action fixing a nonempty finite set $M^{({\Bbb Z}_2)^n}$, which
implies $m\geq n$ (see \cite{ap}). Then we know from \cite{l2, l3}
that the $({\Bbb Z}_2)^n$-action on $M^m$ defines a regular graph $\Gamma_M$ of
valence $m$ with the vertex set $M^{({\Bbb Z}_2)^n}$ and a $({\Bbb Z}_2)^n$-coloring
$\alpha$.
In  the extreme
case where $m=n$ (i.e., $M^n$ is a 2-torus manifold), one knows from \cite{bl, l2, l3} that such a $({\Bbb Z}_2)^n$-colored graph $(\Gamma_M, \alpha)$ is uniquely determined by the
$({\Bbb Z}_2)^n$-action where $\alpha$ is defined as a  map  from the
set $E_{\Gamma_M}$ of all edges of $\Gamma_M$ to all non-trivial
elements of $J_n^{\Bbb R}$, and it satisfies the following
properties:.
\begin{enumerate}
\item[(P1)] for each vertex $v$ of $\Gamma_M$, $\prod_{x\in E_v}\alpha(x)$ is faithful
in ${\Bbb Z}_2[J_n^{\Bbb R}]$, where $E_v$ denotes the set of all edges
adjacent to $v$;
\item[(P2)] for each edge $e$ of $\Gamma_M$, $\alpha(E_u)\equiv
\alpha(E_v) \mod \alpha(e)$ in $J_n^{\Bbb R}$ where $u$ and
$v$ are two endpoints of $e$.
\end{enumerate}
The pair $(\Gamma_M, \alpha)$ is called the {\em $({\Bbb Z}_2)^n$-colored
graph} of the 2-torus manifold $M^n$.
\vskip .1cm

Guillemin and Zara \cite{gz} formulated the results of GKM theory in
terms of a colored graph, and developed the GKM theory
combinatorially. They defined and studied the abstract GKM graphs.
This idea may still be carried out in the mod 2 case.  Following
\cite{l2}, let $\Gamma$ be a finite regular graph of valence $n$
without loops. If there is a map $\alpha$ from the set $E_\Gamma$ of
all edges of $\Gamma$ to all nontrivial elements of $J_n^{\Bbb R}$ satisfying  the properties (P1) and (P2) as above,  then the
pair $(\Gamma, \alpha)$ is called an  {\em abstract $({\Bbb Z}_2)^n$-colored
graph} of $\Gamma$, and $\alpha$ is called a {\em $({\Bbb Z}_2)^n$-coloring} on
$\Gamma$.
\vskip .1cm  Let $(\Gamma, \alpha)$ be an abstract $({\Bbb Z}_2)^n$-colored graph.
Set $$g_{(\Gamma, \alpha)} = \sum_{v\in V_\Gamma} \prod_{x\in
E_v}\alpha(x)$$ which is called the {\em $({\Bbb Z}_2)^n$-coloring polynomial}
of $(\Gamma, \alpha)$.
Obviously, $g_{(\Gamma, \alpha)}$ is a faithful $({\Bbb Z}_2)^n$-polynomial  in ${\Bbb Z}_2[J_n^{\Bbb R}]$. It was shown in \cite[Proposition 2.2]{l2} that for an  abstract
$({\Bbb Z}_2)^n$-colored graph $(\Gamma, \alpha)$, the collection
$\{\alpha(E_v), v\in V_{\Gamma}\}$ is always realizable as the fixed
point data of some 2-torus manifold $M^n$, which implies that the
$({\Bbb Z}_2)^n$-coloring polynomial $g_{(\Gamma, \alpha)}$ of $(\Gamma,
\alpha)$ must belong to the image  $\text{\rm Im}\phi_n$. However, this result does not
tell us whether $(\Gamma, \alpha)$ is the $({\Bbb Z}_2)^n$-colored graph
$(\Gamma_M, \alpha)$ of $M^n$ or not, which is related  to the
following geometric realization problem: {\em under what condition
can $(\Gamma, \alpha)$ become a $({\Bbb Z}_2)^n$-colored graph of some 2-torus
manifold?} Some work for the geometric realization problem has been
studied in details in \cite{bl}.
On the other hand, one has known  from \cite{bl}
or \cite[Section 2]{l3} that each 2-torus manifold $M^n$  determines
a $({\Bbb Z}_2)^n$-colored graph $(\Gamma_M, \alpha)$, and the corresponding
$({\Bbb Z}_2)^n$-coloring polynomial $g_{(\Gamma_M, \alpha)}$ is exactly
$\phi_n(\{M^n\})$.  
This gives another characterization of a faithful $({\Bbb Z}_2)^n$-polynomial $g\in \text{\rm Im}\phi_n$ in terms of $({\Bbb Z}_2)^n$-colored graphs.
\begin{thm} [{\cite[Theorem 4.2]{lt}}]\label{color poly}
A faithful $({\Bbb Z}_2)^n$-polynomial $g$ in ${\Bbb Z}_2[J_n^{\Bbb R}]$ belongs to $\text{\rm Im}\phi_n$ if and only if it is the
$({\Bbb Z}_2)^n$-coloring polynomial of an abstract $({\Bbb Z}_2)^n$-colored graph
$(\Gamma, \alpha)$.
\end{thm}
\begin{rem}
It was shown in~\cite[Propositions 5.2--5.3]{lt} that for a faithful $({\Bbb Z}_2)^n$-polynomial $g\in {\Bbb Z}_2[J_n^{\Bbb R}]$, if $g\in \text{\rm Im}\phi_n$ then  $d(g^*)=0$, and if
$d(g^*)=0$, then $g$ is the $G_n$-coloring polynomial of an abstract
$({\Bbb Z}_2)^n$-colored graph. Then Theorem~\ref{main result} follows from Theorem~\ref{color poly}.
\end{rem}
Now by $\mathcal{G}(({\Bbb Z}_2)^n)$ we denote the set of all abstract
$({\Bbb Z}_2)^n$-colored graphs $(\Gamma, \alpha)$.
Two abstract $({\Bbb Z}_2)^n$-colored graphs $(\Gamma_1, \alpha_1)$ and
$(\Gamma_2, \alpha_2)$ in $\mathcal{G}(({\Bbb Z}_2)^n)$ are said to be {\em
equivalent} if $g_{(\Gamma_1, \alpha_1)}=g_{(\Gamma_2, \alpha_2)}$,
denoted by $(\Gamma_1, \alpha_1)\sim(\Gamma_2, \alpha_2)$.
On the coset $\mathcal{G}(({\Bbb Z}_2)^n)/\sim$, define the addition $+$ as
follows:
$$\{(\Gamma_1, \alpha_1)\}+\{(\Gamma_2, \alpha_2)\}:=\{(\Gamma_1, \alpha_1)\sqcup(\Gamma_2, \alpha_2)\}$$
where $\sqcup$ means the disjoint union.
\begin{cor}[\cite{lt}]\label{graph}
$\mathcal{Z}_n(({\Bbb Z}_2)^n)$ is isomorphic to $\mathcal{G}(({\Bbb Z}_2)^n)/\sim$.
\end{cor}

Davis--Januszkiewicz theory  gives another link between the equivariant topology
and the combinatorics of simple convex polytopes. An $n$-dimensional {\em small cover} $\pi:
M^n\longrightarrow P^n$ is a smooth closed $n$-manifold $M^n$ with a
locally standard $({\Bbb Z}_2)^n$-action such that its orbit space is a simple convex $n$-polytope $P^n$, where a locally
standard $({\Bbb Z}_2)^n$-action on $M^n$ means that this $({\Bbb Z}_2)^n$-action on $M^n$
is locally isomorphic to a faithful representation of $({\Bbb Z}_2)^n$ on
${\Bbb R}^n$.  Each
small cover $\pi: M^n \longrightarrow P^n$ determines a
characteristic function $\lambda$ (here we call it a {\em
$({\Bbb Z}_2)^n$-coloring}) on $P^n$, defined by mapping all facets (i.e.,
$(n-1)$-dimensional faces) of $P^n$ to nontrivial elements of
$J_n^{*{\Bbb R}}$ such that $n$ facets meeting at each vertex are
mapped to $n$ linearly independent elements. A fascinating characteristic for $\pi: M^n\longrightarrow P^n$ is that $M^n$ can be recovered by the pair $(P^n, \lambda)$, so that the algebraic topology of $M^n$ is essentially
consistent with the algebraic combinatorics of $(P^n, \lambda)$.
\vskip .1cm

 Now suppose that $\pi: M^n\longrightarrow P^n$ is a small cover, and
 $\lambda: \mathcal{F}(P^n)\longrightarrow
 J_n^{*{\Bbb R}}$ is its characteristic function, where $\mathcal{F}(P^n)$ consists of all facets of $P^n$.
Given a vertex $v$ of $P^n$, since $P^n$ is simple, there are $n$
facets $F_1, ..., F_n$ in $\mathcal{F}(P^n)$ such that
$v=F_1\cap\cdots\cap F_n$. Then the vertex $v$  determines a monomial
$\prod_{i=1}^n\lambda(F_i)$  of degree $n$ in
${\Bbb Z}_2[J_n^{*{\Bbb R}}]$, whose dual by the pairing
(\ref{pairing}) is faithful in ${\Bbb Z}_2[J_n^{\Bbb R}]$.
Here $\prod_{i=1}^n\lambda(F_i)$ is called the {\em $({\Bbb Z}_2)^n$-coloring
monomial at $v$},  denoted by $\lambda_v$.  Moreover, all vertices
in the vertex set $V_{P^n}$ of $P^n$ via $\lambda$ give a polynomial
  $\sum_{v\in V_{P^n}}\lambda_v$ of degree $n$ in ${\Bbb Z}_2[J_n^{*{\Bbb R}}]$, which is denoted by
 $g_{(P^n, \lambda)}$, and is called  $g_{(P^n, \lambda)}$ the {\em $({\Bbb Z}_2)^n$-coloring polynomial}
 of $(P^n, \lambda)$. On the other hand,  let $(\Gamma_M, \alpha)$ be the $({\Bbb Z}_2)^n$-colored graph of $\pi:
M^n\longrightarrow P^n$, and let $g_{(\Gamma_M, \alpha)}$ be the
$({\Bbb Z}_2)^n$-coloring polynomial of $(\Gamma_M, \alpha)$. One knows from
\cite[Proposition 4.1; Remark 4]{l2} that $\Gamma_M$ is exactly the
1-skeleton of $P^n$, and both $\lambda$ and $\alpha$ determine each
other.  This gives
\begin{prop}[{\cite[Proposition 4.7]{lt}}]\label{small}
$g_{(P^n, \lambda)}$ is the dual polynomial of $g_{(\Gamma_M,
\alpha)}$.
\end{prop}
\begin{rem}\label{p-formula}
The following formulae  for a colored polytope $(P^n, \lambda)$ with $P^n=P_1\times P_2$ were obtained.
\begin{enumerate}
\item[$\bullet$] Product formula (\cite[Proposition 4.10]{lt}). $g_{(P_1\times P_2, \lambda)}=g_{(P_1, \lambda_1)}g_{(P_2,
\lambda_2)};$
\item[$\bullet$]Connected sum formula (\cite[Proposition 4.12]{lt}). $g_{(P_1\sharp_{v_1, v_2} P_2, \lambda)}=g_{(P_1, \lambda_1)}+g_{(P_2,
\lambda_2)}$
\end{enumerate}
where $\lambda_i$ is the restriction to $P_i$ of $\lambda$, and $v_i$ is a vertex of $P_i$.
\end{rem}
Proposition~\ref{small} provides us much insight to  the study on $\mathcal{Z}_n(({\Bbb Z}_2)^n)$. A further question is {\em whether or not can the dual $g^*$ of $g\in \text{\rm Im} \phi_n$ be characterized in terms of $({\Bbb Z}_2)^n$-colored $n$-polytopes?} The positive solution of this question means that the Conjecture $(*)$ holds.

 \subsection{Structure of $\mathcal{Z}_n(({\Bbb Z}_2)^n)$} Now let us look at the structure of $\mathcal{Z}_n(({\Bbb Z}_2)^n)$. One has by Theorem~\ref{main result}  that as linear spaces over ${\Bbb Z}_2$, $\mathcal{Z}_n(({\Bbb Z}_2)^n)$ is   isomorphic to the linear space
  $\mathcal{V}_n$ formed by all faithful $({\Bbb Z}_2)^n$-polynomials $g\in {\Bbb Z}_2[J_n^{\Bbb R}]$  with $d(g^*)=0$. Then, the problem can be further reduced to studying the linear space
   $\mathcal{V}^*_n$ formed by the dual polynomials of those polynomials in $\mathcal{V}_n$. In~\cite[Proposition 6.7]{lt}, L\"u and Tan first showed that $\mathcal{V}^*_n$ is generated by the $({\Bbb Z}_2)^n$-polynomials of products of simplices with $({\Bbb Z}_2)^n$-colorings, and they then showed in the proof of \cite[Theorem 2.5]{lt} that each polynomial of  $\mathcal{V}^*_n$ is exactly the $({\Bbb Z}_2)^n$-coloring polynomial of a $({\Bbb Z}_2)^n$-colored simple convex $n$-polytope. This gives
\begin{thm}[{\cite[Corollary 6.10]{lt}}] \label{relation1}
A faithful $({\Bbb Z}_2)^n$-polynomial $g\in {\Bbb Z}_2[J_n^{\Bbb R}]$ belongs to $\text{\rm Im}\phi_n$ if and only if its dual polynomial
 $g^*$ is the $({\Bbb Z}_2)^n$-coloring polynomial of a $({\Bbb Z}_2)^n$-colored simple convex polytope $(P^n, \lambda)$.
\end{thm}
As a consequence, one has that
\begin{cor}[{\cite[Theorem 2.5]{lt}}]
The Conjecture $(*)$ holds.
\end{cor}
The proof of Theorem~\ref{relation1} in~\cite{lt} also tells us  the basic structure of the
graded noncommutative ring $\mathfrak{M}_*=\sum_{n\geq 1}\mathcal{Z}_n(({\Bbb Z}_2)^n)$, which is stated as follows.
\begin{thm} [{\cite[Theorem 2.6]{lt}}]\label{compute}
$\mathfrak{M}_*$ is generated by the equivariant unoriented bordism classes of all generalized real Bott manifolds.
\end{thm}
\begin{rem}
Generalized real Bott manifolds belong to a class of nicely behaved
small covers, which were introduced  and studied in \cite{cms}. A {\em generalized real Bott tower} of height $n$ is a
sequence of ${\Bbb R}P^{n_i}$-bundles with $n_i\geq 1$:
\[
\begin{CD} B^{\Bbb R}_n@ >{\pi_n}>> B^{\Bbb R}_{n-1}@>{\pi_{n-1}}>> \cdots@>{\pi_2}>>B^{\Bbb R}_{1}@
>{\pi_1}>>B^{\Bbb R}_{0}=\{\text{a point}\}
\end{CD}
\]
where each $\pi_i: B^{\Bbb R}_i\longrightarrow B^{\Bbb R}_{i-1}$ for $i=1, ..., n$ is the projectivization of a Whitney sum of $n_i+1$ real line bundles over $B^{\Bbb R}_i$, and $B^{\Bbb R}_i$ is called an {\em $i$-stage generalized real Bott
manifold} or a {\em generalized real Bott manifold of height $i$}. It is well-known that $B^{\Bbb R}_i$ is a small cover over $\Delta^{n_1}\times\cdots\times \Delta^{n_i}$ where $\Delta^{n_j}$ is an $n_j$-dimensional simplex.
\end{rem}
Theorem~\ref{compute} gives a strong connection between the computation
of equivariant bordism groups $\mathcal{Z}_n(({\Bbb Z}_2)^n)$ and the Davis--Januszkiewicz theory of small covers. As a computational application, one has that
\begin{prop}[\cite{lt}] For $n=3$,
$\dim_{{\Bbb Z}_2}\mathcal{Z}_3(({\Bbb Z}_2)^3)=13$, and for $n=4$, $\dim_{{\Bbb Z}_2}\mathcal{Z}_4(({\Bbb Z}_2)^4)=510$.
\end{prop}

\subsection{A summary and further problems}\label{app-2}

 Together with Theorems~\ref{dks}, \ref{main result}, \ref{color poly} and~\ref{relation1}, it follows that there are some essential relationships among 2-torus manifolds, coloring polynomials, colored simple convex polytopes, colored graphs, which are stated as follows:

 \begin{thm}\label{summary}
 Let $g=\sum_i t_{i,1}\cdots t_{i, n}$ be a faithful $({\Bbb Z}_2)^n$-polynomial in ${\Bbb Z}_2[J_n^{\Bbb R}]$, and $g^*$ be the dual polynomial of $g$. Then the following
 statements are all equivalent.

\begin{enumerate}
\item[(1)] $g\in \text{\rm Im} \phi_n$ $($i.e., there is an $n$-dimensional 2-torus manifold $M^n$ such that $g=\sum_{p\in M^G}[\tau_pM])$;

\item[(2)] $g$ is the $({\Bbb Z}_2)^n$-coloring polynomial of a $({\Bbb Z}_2)^n$-colored graph $(\Gamma, \alpha)$;

\item[(3)] $g=\sum_i t_{i,1}\cdots t_{i, n}$ possesses the property that for any symmetric
polynomial function $f(x_1, ..., x_n)$ over ${\Bbb Z}_2$,
$$\sum_{i}{{f(t_{i,1}, ..., t_{i, n})}\over{t_{i,1}\cdots t_{i, n}}}\in\mathcal{S}(J_n^{\Bbb R});$$

\item[(4)] $d(g^*)=0$;

\item[(5)] $g^*$ is the $({\Bbb Z}_2)^n$-coloring polynomial of a $({\Bbb Z}_2)^n$-colored simple convex polytope $(P^n, \lambda)$.
\end{enumerate}
\end{thm}

 Based upon the above equivalent results, it seems to be interesting to discuss the properties of
 regular graphs and simple convex polytopes. In Theorem~\ref{summary}(2),  $\Gamma$  can actually be chosen as the 1-skeleton of a polytope.
 However, for a $({\Bbb Z}_2)^n$-colored graph $(\Gamma, \alpha)$, we don't know when $\Gamma$ will become the 1-skeleton of a polytope.  Indeed, given a graph, to determine whether it is the 1-skeleton of a polytope or not
 is a quite difficult problem except for the known Steinitz theorem (see \cite{g}). In addition, the product formula in Remark~\ref{p-formula} tells us that  a simple convex polytope with a coloring is indecomposable if its coloring polynomial is indecomposable. These observations lead us  to pose  the following problems:
\begin{enumerate}
\item[$\bullet$] {\em For a $({\Bbb Z}_2)^n$-colored graph $(\Gamma, \alpha)$, under what condition will $\Gamma$ be the 1-skeleton of a polytope?}
\item[$\bullet$] {\em Given a $({\Bbb Z}_2)^n$-colored simple convex polytope $(P^n, \lambda)$, can we give a necessary and sufficient condition that $P^n$ is indecomposable?}
\end{enumerate}


\section{Equivariant unitary bordism  of  unitary toric manifolds}\label{unitary}

So far, there have been many different ways to the study of equivariant unitary bordism ring $\Omega_*^{U, G}$ (see, e.g., \cite{bpr, cf1, tom3, lo, lw, n, n1}), where $G$ is a compact Lie group, and $\Omega_*^{U, G}$ is generated by the equivariant unitary bordism classes of all unitary $G$-manifolds.  However, the explicit ring structure of  $\Omega_*^{U, G}$ is still difficult to calculate and only partial results are known.
Some explicit computations of various equivariant unitary bordism rings were made by Kosniowski~\cite{k}, Landweber~\cite{la}, and Stong~\cite{s2} for $G={\Bbb Z}_p$ and by Kosniowski--Yahia~\cite{ky} and Sinha~\cite{sin} for $G=S^1$. More recently, when $G=T^n$,  Hanke showed in~\cite[Theorem 1]{han} the existence of a pull back square
\[\begin{CD}
  \Omega_*^{U, T^n} @>>> MU_*[e_V^{-1}, Y_{V,d}]\\
  @VVV @VVV\\
 MU_*^{T^n} @>>> MU_*[e_V, e_V^{-1}, Y_{V,d}]
\end{CD}
\]
with all maps being injective, so that the Pontrjagin--Thom map induces a ring isomorphism
$$\Omega_*^{U, T^n}\cong  MU_*^{T^n}\cap MU_*[e_V^{-1}, Y_{V,d}]$$
where  $MU_*^{T^n}$ is the homotopy theoretic equivariant unitary bordism ring (which was defined by tom Dieck in~\cite{tom3}),  $MU_*$ is the ordinary homotopy theoretic unitary bordism ring, and $e_V$ are the Euler classes of nontrivial irreducible complex $T^n$-representations $V$, $Y_{V, d}$ are the classes of degree $2d$ $(2\leq d\leq \infty)$ represented by the $T^n$-bundle $E\otimes V\longrightarrow {\Bbb C}P^{d-1}$ with $E\longrightarrow {\Bbb C}P^{d-1}$ being the hyperplane line bundle (see~\cite{han} for more details).  In his Ph. D's thesis~\cite[Proposition 4.24]{dar}, Darby  showed an refinement result to the case in which the fixed point set is isolated, saying that the commutative diagram
\[\begin{CD}
  \mathcal{Z}_*^{U}(T^n) @>>> \Omega_*^{U, T^n}\\
  @VVV @VVV\\
{\Bbb Z}[e_V^{-1}] @>>> MU_*[e_V^{-1}, Y_{V,d}]
\end{CD}
\]
has all maps injective and is a pullback square, where $\mathcal{Z}_*^{U}(T^n)=\sum_{m\geq 0}\mathcal{Z}_{2m}^{U}(T^n)$ is the subring of $\Omega_*^{U, T^n}$ that consists of all classes that can be represented by unitary $T^n$-manifolds with finite fixed point set. Furthermore, Darby obtained  in~\cite[Corollary 6.20]{dar} the following monomorphism from his refinement result:
$$\varphi_n: \mathcal{Z}_{2n}^{U}(T^n)\longrightarrow \Lambda_{\Bbb Z}^n(J^{\Bbb C}_n)$$
where $J^{\Bbb C}_n$ denotes the set of  irreducible $T^n$-representations so $J^{\Bbb C}_n$ can be identified with $ \text{\rm Hom}(T^n, S^1)$, and $\Lambda_{\Bbb Z}(J^{\Bbb C}_n)$ is the free exterior algebra on $J^{\Bbb C}_n$ over ${\Bbb Z}$ with the graded structure $\Lambda_{\Bbb Z}(J^{\Bbb C}_n)=\bigoplus_m\Lambda^m_{\Bbb Z}(J^{\Bbb C}_n)$ having the property $\Lambda^k_{\Bbb Z}(J^{\Bbb C}_n)\wedge\Lambda^l_{\Bbb Z}(J^{\Bbb C}_n)\subset \Lambda^{k+l}_{\Bbb Z}(J^{\Bbb C}_n)$. The monomorphism $\varphi_n: \mathcal{Z}_{2n}^{U}(T^n)\longrightarrow \Lambda^n_{\Bbb Z}(J^{\Bbb C}_n)$ is an analogue of $\phi_*: \mathcal{Z}_*(({\Bbb Z}_2)^n)\longrightarrow {\Bbb Z}_2[J_n^{\Bbb R}]$ as stated in Section~\ref{2-torus}.
Moreover, Darby carried out his study on $\mathcal{Z}_{2n}^{U}(T^n)$ (i.e., the equivariant unitary bordism of all unitary toric $2n$-manifolds) by capturing the ideas developed in the setting of 2-torus manifolds, which will be introduced in the next subsection.

 \subsection{Structure of $\Xi_*$, faithful exterior polynomials and torus graphs}
An exterior polynomial $g$ in $\Lambda_{\Bbb Z}^n(J^{\Bbb C}_n)$ is said to be {\em faithful} if the indeterminates
from each monomial of $g$ form a basis of $J^{\Bbb C}_n$. Similarly to the definition of $\Lambda_{\Bbb Z}(J^{\Bbb C}_n)$, define $\Lambda_{\Bbb Z}(J_n^{*{\Bbb C}})$ to be the free exterior algebra on $J_n^{*{\Bbb C}}$ over ${\Bbb Z}$, where $J_n^{*{\Bbb C}}=\text{\rm Hom}(S^1, T^n)$. Since both $J_n^{{\Bbb C}}$ and $J_n^{*{\Bbb C}}$ are isomorphic to ${\Bbb Z}^n$ and are dual by the following pair
\begin{equation*}\label{pairing1}
\langle\cdot, \cdot\rangle: J_n^{*{\Bbb C}}\times J_n^{{\Bbb C}}\longrightarrow \text{\rm Hom}(S^1, S^1)\cong {\Bbb Z}
\end{equation*}
defined by $\langle \xi, \rho\rangle=\rho\circ \xi$, for each faithful exterior polynomial $g\in \Lambda_{\Bbb Z}^n(J^{\Bbb C}_n)$ one may obtain a {\em dual polynomial} $g^*\in \Lambda_{\Bbb Z}^n(J_n^{*{\Bbb C}})$ by considering the dual basis in $J_n^{*{\Bbb C}}$ of the basis in $J^{\Bbb C}_n$ produced by each monomial of $g$. Clearly, the monomorphism $\varphi_n$ maps each nonzero class $\beta$ of $\mathcal{Z}_{2n}^{U}(T^n)$ to an faithful exterior polynomial $\varphi_n(\beta)$ in $\Lambda_{\Bbb Z}^n(J^{\Bbb C}_n)$.

\vskip .1cm
Similarly, a differential operator $d$ on $\Lambda_{\Bbb Z}(J_n^{*{\Bbb C}})$ may be defined as follows: for each monomial $s_1\wedge\cdots\wedge s_k\in \Lambda^k_{\Bbb Z}(J_n^{*{\Bbb C}})$ with $k\geq 1$
$$d_k(s_1\wedge\cdots\wedge s_k)=\begin{cases}
\sum_{i=1}^k(-1)^{i+1}s_1\wedge\cdots s_{i-1}\wedge \widehat{s}_i\wedge s_{i+1}\wedge\cdots\wedge s_k & \text{\rm if } k>1\\
1 & \text{\rm if } k=1
\end{cases}
$$
and $d_0(1)=0$.
\vskip .1cm
Recall that $\Xi_*=\bigoplus_{n\geq 0}\mathcal{Z}_{2n}^{U}(T^n)$ is the graded noncommutative ring generated by the equivariant unitary bordism classes of all unitary toric manifolds. Let $K_n$ denote the  abelian group of all faithful exterior polynomials
$g\in \Lambda_{\Bbb Z}^n(J^{\Bbb C}_n)$ such that $d(g^*)=0$. Then $K_*=\bigoplus_{n\geq 0}K_n$ forms a graded noncommutative subring of $\Lambda_{\Bbb Z}(J^{\Bbb C}_n)$. Darby showed that
\begin{thm}[{\cite[Theorem 8.5]{dar}}] Both $\Xi_*$ and $K_*$ are isomorphic. Furthermore, a faithful exterior polynomial $g$ in $\Lambda_{\Bbb Z}^n(J^{\Bbb C}_n)$ belongs to $\text{\rm Im}\varphi_n$ if and only if $d(g^*)=0$.
\end{thm}

A {\em torus graph} is a pair $(\Gamma, \alpha)$ consisting of an $n$-valent regular graph $\Gamma$ with
a torus axial function $\alpha: E_\Gamma\longrightarrow J^{\Bbb C}_n$ subject to the following properties:
\begin{enumerate}
\item $\alpha(\overline{e})=\pm\alpha(e)$;
\item for each vertex $v$, $\alpha(E_v)$ forms a basis of $J^{\Bbb C}_n$;
\item for each edge $e$, $\alpha(E_{i(e)})\equiv \alpha(E_{t(e)})\mod \alpha(e)$
\end{enumerate}
where $E_\Gamma$ denotes the set of oriented edges of $\Gamma$, that is, each edge appears twice in
$E_\Gamma$ with opposite orientations, $i(e)$ and $t(e)$ denote the initial and terminal vertices of an edge $e\in E_\Gamma$, respectively,  $\overline{e}$ denotes the edge $e$ with its opposite orientation, and for each vertex $v$,
$E_v=\{e\in E_\Gamma| i(e)=v\}$. Note that a torus axial function is different from the axial function for
GKM-graphs as defined in \cite{gz1}, which requires that $\alpha(\overline{e})=-\alpha(e)$ as well as that
the elements of $\alpha(E_v)$ are pairwise linearly independent. As shown in~\cite{mmp}, a torus graph is not a GKM graph in general, each torus manifold $M$ determines a torus graph $(\Gamma_M, \alpha)$, and all torus graphs are orientable, where an {\em orientation} of a torus graph $(\Gamma, \alpha)$ is an assignment $\sigma: V_\Gamma\longrightarrow\{\pm1\}$
satisfying $\sigma(i(e))\alpha(e)=-\sigma(i(\overline{e})\alpha(\overline{e})$  for every $e\in E_\Gamma$, and a {\em torus manifold} of dimension $2n$ is a smooth closed $2n$-manifold with an effective $T^n$-action fixing a nonempty set (so a unitary toric manifold is a special torus manifold). Darby showed in~\cite[Proposition 6.11]{dar} that the torus graph of a unitary toric manifold is orientable, and he gave the definition of the torus polynomial
 $g_{(\Gamma, \alpha, \sigma)}$ of an oriented torus graph $(\Gamma, \alpha, \sigma)$ (see~\cite[Definition 6.17]{dar}).  Furthermore, Darby characterized the torus polynomials of  oriented torus graphs in terms of the vanishing of the differential $d$ on the dual polynomials.

\begin{thm}[{\cite[Theorem 6.26]{dar}}]
Let $g\in \Lambda_{\Bbb Z}^n(J^{\Bbb C}_n)$ be a faithful polynomial. Then $g$ is the torus polynomial of an oriented torus graph if and only if $d(g^*)=0$.
\end{thm}

A {\em quasitoric manifold} is an even-dimensional smooth closed manifold
$M^{2n}$ equipped  with a locally standard smooth $T^n$-action such that the orbit space is a simple
$n$-polytope $P$. Like small covers, each quasitoric manifold $\pi: M^{2n}\longrightarrow P^n$ determines a characteristic map $\lambda$ (also called a ${\Bbb Z}^n$-coloring here) on $P^n$, which sends each facet of $P^n$ onto non-trivial elements of $J^{*{\Bbb C}}_n$, unique up to sign, such that the $n$ facets of $P^n$ meeting at a vertex are mapped to a basis of $J^{*{\Bbb C}}_n$, and $M^{2n}$ can be recovered from the combinatorial data $(P^n, \lambda)$.   As shown in~\cite{br}, each quasitoric manifold $\pi: M^{2n}\longrightarrow P^n$ with an omniorientation is a unitary toric manifold, where an {\em omniorientation} consists of a choice of orientation for $M^{2n}$ and for every  submanifold $\pi^{-1}(F)$, $F\in \mathcal{F}(P^n)$ (the set of all facets of $P^n$). Each omnioriented quasitoric manifold still determines a pair $(P^n, \lambda)$, which is called the {\em quasitoric pair}. Then the following result means that the study on omnioriented quasitoric manifolds can be reduced to the study on quasitoric pairs.
\begin{thm}[{\cite[Corollary 5.28]{dar}, also see~\cite[Theorem 5.10]{bpr}}]
There is a bijection between the set of quasitoric pairs and the set
of omnioriented quasitoric manifolds.
\end{thm}

To consider whether  each class of $\mathcal{Z}_{2n}^{U}(T^n)$ is represented by an omnioriented quasitoric manifold, Darby in~\cite{dar} introduced a graded noncommutative ring $\mathcal{Q}_*$, which is  generated by all quasitoric pairs with the addition and the multiplication given by the disjoint union and the cartesian product, respectively. Then he studied the homomorphism of non-commutative graded rings
$$\mathcal{M}: \mathcal{Q}_*\longrightarrow \Xi_*$$
by constructing the omnioriented quasitoric manifold associated to a quasitoric pair. This homomorphism is not a monomorphism, but if it is surjective, then one can obtain that each class of $\mathcal{Z}_{2n}^{U}(T^n)$ is represented by an omnioriented quasitoric manifold. Darby made a significant advance to show that when $n=1,2$, $\mathcal{M}$ is surjective (see~\cite[Corollaries 8.8 and 8.10]{dar}).
\vskip .1cm
 Darby also extended the connected
sum construction of quasitoric pairs which allows for a more general notion of the
equivariant connected sum of omnioriented quasitoric manifolds (\cite[\S 7.3]{dar}), and obtained the connected sum formula and the product formula for quasitoric pairs (\cite[Lemmas 7.5 and 7.11]{dar}). In addition, for an omnioriented quasitoric manifold $\pi: M^{2n}\longrightarrow P^n$,
he also showed that the polynomial of the quasitoric pair $(P^n, \lambda)$ is the dual of the torus polynomial of the associated oriented torus graph $(\Gamma_M, \alpha)$ (\cite[Formula (7.5)]{dar}).

\subsection{Equivariant Chern numbers and the number of fixed points for unitary torus manifolds}
In~\cite{ggk}, Guillemin,  Ginzburg and  Karshon showed that the equivariant
unitary bordism class of a unitary $T^n$-manifold with
isolated fixed points is completely determined by its equivariant
Chern numbers. In~\cite{lt1}, L\"u and Tan gave a refinement of their result in the setting of unitary toric manifolds.
\begin{thm}[{\cite[Theorem 1.1]{lt1}}]\label{bounds}
Let $\beta=\{M\}$ be a class in $\mathcal{Z}_{2n}^U(T^n)$. Then~$\beta=0$
 if and only if the equivariant Chern
numbers $\langle
(c_1^{T^n})^i(c_2^{T^n})^j, [M]\rangle=0$ for all $i, j\in {\Bbb
N}$, where $[M]$ is the fundamental class of $M$ with respect to the given orientation.
\end{thm}

In~\cite{ck}, Kosniowski studied unitary~$S^1$-manifolds and  when
the fixed points are isolated, he proposed the following conjecture.
\vskip .1cm
\noindent {\bf Conjecture} (Kosniowski). {\em Suppose that~$M^{2n}$ is a unitary $S^1$-manifold with isolated fixed points. If~$M^{2n}$ does not bound equivariantly then the number of fixed points is greater than $f(n)$, where $f(n)$ is some linear function.}

\begin{rem}
As was noted by Kosniowski in \cite{ck},  the most likely function is $f(n)={n\over 2}$, so the number of fixed points
of $M^{2n}$ is at least~$[{n\over 2}]+1$.
\end{rem}

With respect to this conjecture, recently some related works have been done (see~\cite{ckp,ll, pt}).
For example, Li and
Liu showed  in~\cite{ll} that if~$M^{2mn}$ is an almost complex manifold and there
exists a partition~$\lambda=(\lambda_1,\dots,\lambda_r)$ of weight
$m$ such that the corresponding Chern number~$\langle(c_{\lambda_1}\dots
c_{\lambda_r})^n, [M]\rangle$ is nonzero, then  any~$S^1$-action on~$M$ must have at least~$n+1$ fixed points.
In the case of the unitary torus manifolds,   one can apply Theorem~\ref{bounds} to obtain the following result, which further provides supporting evidence to the Kosniowski conjecture.

\begin{thm}[{\cite[Theorem 1.2]{lt1}}]\label{number}
Suppose that~$M^{2n}$ is a $(2n)$-dimensional unitary toric manifold. If~$M^{2n}$ does not bound equivariantly, then
the number of fixed points is at least~$\lceil{n\over2}\rceil+1$, where $\lceil{n\over2}\rceil$ denotes the minimal integer no less than ${n\over 2}$.
\end{thm}

\subsection{Problems and conjectures} We would like to conclude this section with the following problems and conjectures:
\begin{prob}
Let $g=\sum_i t_{i,1}\wedge\cdots\wedge t_{i, n}$ be a faithful exterior polynomial
in $\Lambda_{\Bbb Z}^n[J_n^{\Bbb C}]$. Is the following fact true?
\vskip .1cm
\noindent{\bf Fact.}   $d(g^*)=0$ if and only if for all symmetric
polynomial functions $f(x_1,...,x_n)$ over ${\Bbb Z}$,
\begin{equation*}\label{integral}
\sum_{i}{{f(t_{i,1}, ..., t_{i, n})}\over{t_{i,1}\cdots t_{i,
n}}}\in\mathcal{S}(J_n^{\Bbb C})\end{equation*} when $t_{i,1}\cdots
t_{i, n}$ and $f(t_{i,1}, ..., t_{i, n})$ are regarded as
polynomials in $\mathcal{S}(J_n^{\Bbb C})$, where $\mathcal{S}(J_n^{\Bbb C})$ is the symmetric algebra on $J_n^{\Bbb C}$ over ${\Bbb Z}$.
\end{prob}

\begin{conj}[{\cite[Conjecture 8.13]{dar}}]
Each class of $\mathcal{Z}_{2n}^{U}(T^n)$ is represented by an omnioriented quasitoric manifold. Equivalently, the homomorphism $\mathcal{M}: \mathcal{Q}_*\longrightarrow \Xi_*$ is surjective.
\end{conj}

\begin{conj}[{\cite[Conjecture 3.1]{lt1}}]
The number $\lceil{n\over2}\rceil+1$ is the best possible lower bound of the number of fixed
points for nonbounding unitary toric manifolds of dimension $2n$.
\end{conj}

\section{Relation between $\mathcal{Z}_n(({\Bbb Z}_2)^n)$ and $\mathcal{Z}_{2n}^{U}(T^n)$}\label{rela}

Milnor's work tells us in \cite{m} (see also \cite{s1}) that there is a  homomorphism
$
F_*: \Omega^U_*\longrightarrow \mathfrak{N}^2_*
$
where  $\mathfrak{N}^2_*=\{\alpha^2| \alpha\in \mathfrak{N}_*\}$. This actually implies that there is a covering homomorphism
$
H_n: \Omega^U_{2n}\longrightarrow \mathfrak{N}_{n}
$
which is induced by $\theta_n\circ F_n$, where $\theta_n:\mathfrak{N}^2_{n}\longrightarrow\mathfrak{N}_n$ is defined by mapping $\alpha^2\longmapsto \alpha$. This observation leads to the following natural question:
{\em Is there a homomorphism $\widetilde{H}_n: \mathcal{Z}_{2n}^{U}(T^n)\longrightarrow \mathcal{Z}_n(({\Bbb Z}_2)^n)$   such that $\widetilde{H}_n$ is onto?}
\vskip .1cm
In~\cite{lt2}, L\"u and Tan discussed this question. The homomorphism $\widetilde{H}_n$ is defined as follows: First,  a class $\{M^{2n}\}$ in $\mathcal{Z}_{2n}^{U}(T^n))$ gives an  oriented torus graph $(\Gamma_M, \alpha)$ of $M^{2n}$. Next one may obtain a $({\Bbb Z}_2)^n$-colored graph $(\Gamma, \overline{\alpha})$ from  $(\Gamma_M, \alpha)$ such that
$\Gamma=\Gamma_M$ and $\overline{\alpha}$ is the mod 2 reduction of $\alpha$, and then the coloring polynomial of $(\Gamma, \overline{\alpha})$ determines a class of $\mathcal{Z}_n(({\Bbb Z}_2)^n)$ as desired. In particular, if $M^{2n}$ is an omnioriented quasitoric manifold over a simple convex polytope $P$, then $\widetilde{H}_n$ exactly maps $\{M^{2n}\}$ into the class of the fixed point set (as a small cover over $P$) of the natural conjugation involution on $M^{2n}$ (see~\cite[Corollary 1.9]{dj}).  Furthermore, it was showed by using Atiyah--Bott--Berline--Vergne localization theorem and a classical result of Stong that $\widetilde{H}_n$ is well-defined. With a very technical method, L\"u and Tan showed that
 \begin{thm}[\cite{lt2}]\label{onto}
The homomorphism $\widetilde{H}_n: \mathcal{Z}_{2n}^{U}(T^n)\longrightarrow \mathcal{Z}_n(({\Bbb Z}_2)^n)$ is surjective.
 \end{thm}
 \begin{rem}
 Theorem~\ref{onto} gives an answer to the lifting problem from small covers to quasitoric manifolds in the sense of equivariant bordism, where the lifting problem is explained as follows: Given a $({\Bbb Z}_2)^n$-coloring $\lambda:
 \mathcal{F}(P^n)\longrightarrow J_n^{*{\Bbb R}}\cong ({\Bbb Z}_2)^n$ on a simple convex $n$-polytope $P^n$, {\em does there exist a
 ${\Bbb Z}^n$-coloring $\widetilde{\lambda}: \mathcal{F}(P^n)\longrightarrow J_n^{*{\Bbb C}}\cong {\Bbb Z}^n$  such that the following diagram commutes?}
 \[ \xymatrix{
&J_n^{*{\Bbb C}} \ar[d]^{\mod 2}\\
\mathcal{F}(P^n) \ar[ur]^{\widetilde{\lambda}}\ar[r]_{\lambda} & J_n^{*{\Bbb R}} }\]
where $\mathcal{F}(P^n)$ denotes the set of all facets of $P^n$. This problem was posed by the author of this paper at the conference on toric topology held in Osaka in November 2011\footnote[2]{http://www.sci.osaka-cu.ac.jp/~masuda/toric/torictopology$2011_-$osaka.html}. The problem is still open except for the cases of $n\leq 3$ and $m-n\leq 3$, where $m$ is the number of all facets of $P^n$ (see~\cite{cp}).
 \end{rem}


\begin{thebibliography}{99}

\bibitem{a}  Alexander, J. C. The bordism ring of manifolds with involution. Proc. Amer. Math. Soc. {\bf 31} (1972),  536--542.
\bibitem{ap} C. Allday and V. Puppe, {\em Cohomological Methods in Transformation Groups}, Cambridge Studies in
Advanced Mathematics, {\bf 32}, Cambridge University Press, 1993.
\bibitem{bl} Z. Q. Bao and Z. L\"u, {\em Manifolds associated with $(\Z_2)^n$-colored regular graphs},  Forum Math {\bf 24} (2012), 121--149.
\bibitem{bgh} D. Biss, V. Guillemin and T. S. Holm,  {\em The mod 2
cohomology of fixed point sets of anti-symplectic involutions}, Adv. Math.  {\bf 185} (2004), 370--399.
\bibitem{bp} V. M. Buchstaber and T.E. Panov, {\em Torus actions and their
applications in topology and combinatorics}, University Lecture
Series, 24. American Mathematical Society, Providence, RI, 2002.

\bibitem{bpr} V. M. Buchstaber, T.E. Panov and N. Ray, {\em Toric Genera}, Internat.
Math. Res. Notices {\bf 2010}, No. 16, 3207--3262.

\bibitem{br} V. M. Buchstaber and N. Ray, {\em Toric manifolds and complex cobordisms},
Uspekhi Mat. Nauk {\bf 53} (1998), 139--140. In Russian; translated in Russ. Math. Surv.
{\bf 53} (1998), 371--373.

\bibitem{ckp} H. W. Cho, J. H. Kim and H. C. Park, {\em On the conjecture of Kosniowski},  Asian J. Math. {\bf 16} (2012), no. 2, 271--278.
\bibitem{cms} S. Y. Choi, M. Masuda and D. Y. Suh, {\em Topological classification of generalized Bott manifolds},
Trans. Amer. Math. Soc. {\bf 362} (2010),  1097--1112.
\bibitem{cp} S. Y. Choi, and H. C. Park, {\em Wedge operations and torus symmetries},  arXiv:1305.0136
\bibitem{co} P.E. Conner, Differentiable Periodic Maps, 2nd Edition, {\bf 738}, Springer--Verlag, 1979.
\bibitem{cf} P.E. Conner and E.E. Floyd, {\em Differentiable  periodic maps}, Ergebnisse Math. Grenzgebiete, N. F., Bd. {\bf 33},  Springer-Verlag, Berlin, 1964.
\bibitem{cf1} P.E. Conner and E.E. Floyd, {\em Periodic maps which preserve a complex structure}, Bull. Amer. Math. Soc. {\bf 70} (1964), 574--579.

\bibitem{dar} A. Darby, {\em Quesitoric manifolds in equivariant complex bordism}, Ph.D. thesis, The University of Manchester, 2013.
\bibitem{dj} M. Davis and T. Januszkiewicz, {\em Convex polytopes, Coxeter orbifolds and torus actions}, Duke Math. J.
 {\bf 61} (1991), 417-451.
\bibitem{tom3} T. tom Dieck, {\em Bordism of $G$-manifolds and integrality theorems}, Topology {\bf 9} (1970), 345--358.
\bibitem{tom1} T. tom Dieck, {\em Characteristic numbers of G-manifolds I},   Invent. Math. {\bf 13} (1971), 213--224.
\bibitem{ggk} V. Guillemin, V. Ginzburg and Y. Karshon, {\em Moment maps, cobordisms, and Hamiltonian group actions}, in `Mathematical Surveys and Monographs', {\bf 98}, American Mathematical Society, Providence, RI, 2002.
 \bibitem{gkm} M. Goresky, R. Kottwitz, R. MacPherson, \emph{Equivariant cohomology, Koszul duality, and the localization theorem}, Invent. Math. \textbf{131} (1998), 25--83.
 \bibitem{g} B. Gr\"unbaum, {\em Convex Polytopes}, Second Edition, Graduate Texts in Mathematics
{\bf 221}, Springer--Verlag, 2003.
  \bibitem{gz} V. Guillemin and C. Zara, {\em 1-Skeleta, Betti numbers, and equivariant cohomology},
 Duke Math. J. {\bf 107} (2001), 283--349.
\bibitem{gz1} V. Guillemin and C. Zara,  {\em Equivariant de Rham Theory and
Graphs},   Asian J. Math. {\bf 3} (1999), 49--76.
 \bibitem{han} B. Hanke, {\em Geometric versus homotopy theoretic equivariant bordism}, Math. Ann. {\bf 332} (2005), 677--696.
\bibitem{k} C. Kosniowski, {\em Generators of the ${\Bbb Z}/p$ bordism ring: serendipity},  Math. Z. {\bf 149} (1976), 121--130.
    \bibitem{ck} C. Kosniowski, {\em Some formulae and conjectures associated with circle
actions}, Topology Symposium, Siegen 1979 (Proc. Sympos., Univ.
Siegen, Siegen, 1979), pp331--339, Lecture Notes in Math., {\bf 788},
Springer, Berlin, 1980.

\bibitem{ks} C. Kosniowski and R.E. Stong, {\em $({\Bbb Z}_2)^k$-actions and characteristic numbers}, Indiana Univ. Math. J. {\bf 28} (1979), 723--743.
\bibitem{ky} C. Kosniowski and M. Yahia,  {\em Unitary bordism of circle actions},  Proceedings of the
Edinburgh Mathematical Society {\bf 26} (1983),  97--105.

\bibitem{la} P.S. Landweber, {\em  Equivariant bordism and cyclic groups},  Proc. Amer.
Math. Soc. {\bf 31} (1972),  564--570.
\bibitem{ll} P. Li and K. F. Liu, {\em Some remarks on circle action on
manifolds}, Mathematical Research Letters {\bf 18} (2011), 437--446.
\bibitem{lo} P. L\"offler,  {\em Characteristic Numbers of Unitary Torus-Manifolds},  Bull. Amer. Math. Soc. {\bf 79} (1973), 1262--1263.
\bibitem{lo1}  P. L\"offler, {\em Bordismengruppen Unit\"arer Torusmannigfaltigkeiten},  Manuscripta Mathematica {\bf 12} (1974), 307--327.
\bibitem{l2} Z. L\"u, {\em Graphs of 2-torus actions},   Toric topology, 261--272,
Contemp. Math., {\bf 460}, Amer. Math. Soc., Providence, RI, 2008.
 \bibitem{l1} Z. L\"u, {\em 2-torus manifolds, cobordism and small covers}, Pacific J. Math. {\bf 241} (2009), 285--308.
   \bibitem{l3} Z. L\"u, {\em Graphs and $({\Bbb Z}_2)^k$-actions}, arXiv:
math.AT/0508643.
     \bibitem{lm} Z. L\"u and M. Masuda, {\em Equivariant classification of 2-torus manifolds},  Colloq. Math. {\bf 115}
 (2009), 171--188.
  \bibitem{lt1} Z. L\"u and Q. B. Tan, {\em Equivariant Chern numbers and the number of fixed points for unitary torus manifolds},  Math. Res. Lett. {\bf 18} (2011), no. 6, 1319--1325.
 \bibitem{lt} Z. L\"u and Q. B. Tan, {\em Small covers and the equivariant bordism classification of 2-torus manifolds}, Int. Math. Res. Notices (First published online: September 3, 2013), doi: 10.1093/imrn/rnt183.
    \bibitem{lt2} Z. L\"u and Q. B. Tan, {\em The relation between  equivariant bordism groups of 2-torus manifolds and unitary toric manifolds}, preprint.
    \bibitem{ly}  Z. L\"u and L. Yu, {\em Topological types of 3-dimensional small covers},  Forum Math. {\bf 23}
  (2011),  245--284.
  \bibitem{lw} Z. L\"u and W. Wang, {\em Equivariant cohomology Chern numbers determine equivariant unitary bordism}, preprint.
     \bibitem{ma} M. Masuda, {\em Unitary toric manifolds, multi-fans and equivariant index}, Tohoku Math. J.
{\bf 51} (1999), 237--265.
\bibitem{mmp} H. Maeda, M. Masuda and T. Panov, {\em Torus Graphs and Simplicial
Posets},  Adv. in Math. {\bf 212} (2007), 458--483.
\bibitem{n} S. P. Novikov, {\em Methods of algebraic topology from the point of view of cobordism theory}, Izvestiya Akademii Nauk SSSR. Seriya Matematicheskaya {\bf 32} (1968), 1245--1263 (Russian). English translation in Mathematics of the USSR--Izvestiya {\bf 2} (1968), 1193--1211.
    \bibitem{n1} S. P. Novikov, {\em Adams operations and fixed points}, Izvestiya Akademii Nauk SSSR. Seriya Matematicheskaya {\bf 31} (1967), 855--951 (Russian). English translation in Mathematics of the USSR--Izvestiya {\bf 1} (1967), 827--913.
\bibitem{pt} A. Pelayo and S. Tolman, {\em Fixed points of symplectic
periodic flows}, Ergodic Theory and
Dynamical Systems {\bf 31} (2011), 1237--1247.
\bibitem{m} J.W. Milnor, {\em On the Stiefel--Whitney numbers of complex manifolds and of spin manifolds}, Topology, {\bf 3} (1965), 223--230.

\bibitem{sinh} Dev P. Sinha, {\em Real equivariant bordism and stable transversality obstructions for ${\Bbb Z}/2$},  Proc. Amer. Math. Soc. {\bf 130} (2002), no. 1, 271--281.

\bibitem{sin} Dev P. Sinha, {\em Computations of complex equivariant bordism rings}, Amer. J. Math. {\bf 123} (2001), no. 4, 577--605.

    \bibitem{s1} R.E. Stong, {\em Notes on cobordism theory}, Princeton University Press, 1968.
    \bibitem{s} R.E. Stong, {\em Equivariant bordism and $({\Bbb Z}_2)^k$-actions}, Duke Math. J. {\bf 37} (1970), 779--785.

    \bibitem{s2} R.E. Stong, {\em Complex and Oriented Equivariant Bordism},  In Topology of Manifolds (Proceedings
of The University of Georgia Institute, Athens, GA, 1969), 291--316. Chicago, IL: Markham, 1970.
  \end{thebibliography}
\end{document}